\documentclass[a4paper,10pt,leqno]{amsart}

\title[Non-trivial self-concordances]{Non-trivial self-concordances \\and a recent conjecture by Botvinnik}
\author{Wolfgang Steimle}
\address{Universit\"at Bonn\\
               Mathematisches Institut\\
               Endenicher Allee~60,
               D-53115 Bonn, Germany}
\email{steimle@math.uni-bonn.de}
\date{\today}

\usepackage[active]{srcltx}
\usepackage{pdfsync}
\usepackage{hyperref}
\usepackage{enumerate,amssymb, amsmath}
\usepackage{graphicx}
\usepackage[arrow,curve,matrix,tips,2cell]{xy}
  \SelectTips{eu}{10} \UseTips
  \UseAllTwocells

\DeclareMathAlphabet{\matheurm}{U}{eur}{m}{n}
\DeclareMathOperator{\id}{id}

  \newcommand{\IZ}{\mathbb{Z}}

\theoremstyle{plain}
\newtheorem{theorem}{Theorem}[section]

\newtheorem{corollary}[theorem]{Corollary}

\newtheorem{conjecture}[theorem]{Conjecture}

\theoremstyle{definition}

\newtheorem{remark}[theorem]{Remark}

\theoremstyle{remark}

\newcommand{\arincl}{\ar@{^{(}->}}
\newcommand{\arinclinv}{\ar@{_{(}->}}

\newcommand{\xycomsquare}[8]                % commutative square (xy-Version)
{$$\xymatrix{#1 \ar[r]^-{#2} \ar[d]^{#4} &
    #3 \ar[d]^{#5}  \\
    #6\ar[r]^-{#7} & #8 }$$}

\newcounter{commentnr}

\begin{document}

\begin{abstract}
The goal of this note is to construct, on many manifolds, non-trivial concordances from the identity to itself. This produces counterexamples to a recent conjecture by Botvinnik.
\end{abstract}

\maketitle

\section{Statement of the results}

Recall that a (smooth) concordance on a smooth manifold $M$ is a diffeomorphism of $M\times I$ which is the identity in a neighborhood of $M\times 0\cup \partial M\times I$. For a concor\-dance $H$ of $M$, denote by $e(H)$ the induced diffeomorphism on $M\times 1$. We say that $H$ is \emph{trivial} if it is isotopic to the identity, via an isotopy that fixes a neighborhood of $M\times 0\cup \partial M\times I$.

Concordances can be described, in a stable range, by algebraic $K$-theory. In this note we explain how to use this relationship to prove:

\begin{theorem}\label{thm:main_result}
For $n\geq 9$ there exists a non-trivial concordance $H$ of $S^1\times D^{n-1}$ such that $e(H)=\id$.
\end{theorem}

In fact, more generally we have:

\begin{theorem}\label{thm:main_general_result}
On any smooth compact orientable manifold $M$ of dimension $n\geq 9$ such that $\pi_1(M)=\IZ$, there is a non-trivial concordance $H$ such that $e(H)=\id$. 
\end{theorem}

Recently, in his impressive preprint \cite{Botvinnik(2012)}, Botvinnik has proposed the following ``Topological Conjecture''.

\begin{conjecture}
Let $M$ be a closed manifold equipped with a metric $g$ of positive scalar curvature. If $H$ is a non-trivial concordance on $M$, then $g$ and $e(H)^*g$ are non-isotopic as metrics of positive scalar curvature.
\end{conjecture}

\begin{corollary}
The Topological Conjecture does not hold.
\end{corollary}

\section{Proof of Theorem \ref{thm:main_result}}

Write $M=S^1\times D^{n-1}$ and denote by $C(M)$ the group of concordances modulo isotopy. Let $h\in C(M)$. Shrinking the interval $I=[0,1]$ to $[0, \frac12]$, we may consider $h$ as a self-diffeomorphism of $M\times [0,\frac12]$ which we may extend again to the whole of $M\times I$ by ``flipping'':
\[H\vert_{M\times [\frac12, 1]}:= i\circ h\circ i\]
where $i$ is induced by reflection of $I$ at $\frac12$. Clearly $e(H)=\id$.

\begin{remark}
A further analysis using the Hatcher spectral sequence and surgery theory shows that on the manifold $M$ any concordance from the identity to itself is isotopic to one of this form.
\end{remark}

Note that in the abelian group $C(M)$, we have $H=h+ \tau h$ where $\tau$ denotes the canonical involution on $C(M)$. To prove Theorem \ref{thm:main_result} we therefore need to show that $\tau\neq -\id$.

By the stable parametrized $h$-cobordism theorem \cite{Waldhausen-Jahren-Rognes(2008)} there is a short exact sequence
\[0\to \pi_2^s(M\amalg \{*\}) \to \pi_2 A(M) \xrightarrow{\pi} C(M) \to 0\]
provided $\dim(M) \geq 9$ (this ``stable range'' is due to \cite{Igusa(1988)}). Here $A(M)$ denotes Waldhausen's $K$-theory of spaces \cite{Waldhausen(1985)}; by homotopy invariance we have $\pi_2 A(M)\cong \pi_2 A(S^1)$. By \cite{Vogell(1985)} the functor $A(-)$ carries a canonical involution $T$ so that the map $\pi$ is equivariant up to the sign $(-1)^n$. (Here we use that $M$ is parallelizable.)

By the fundamental theorem \cite{Huettemann-Klein-etc(2001)} 
\begin{equation}\label{eq:fund_theorem}
\pi_2 A(S^1) \cong \pi_2 A(*) \oplus \pi_1 A(*) \oplus \pi_2 NA(*) \oplus \pi_2 NA(*);
\end{equation}
the involution $T$ interchanges the two copies of $\pi_2 NA(*)$ by \cite{Huettemann-Klein-etc(2002)}. Moreover $\pi_k A(*)\cong \pi_k^s$ for $k\leq 2$ in a way that 
\[\pi_2 A(*) \oplus \pi_1 A(*) = \pi_2^s(S^1\amalg \{*\})\]
as subgroups of $\pi_2 A(S^1)$. So 
\[C(M)\cong \pi_2 NA(*)\oplus \pi_2 NA(*)\]
and $\tau$ acts, up to sign, by interchanging the summands. In particular $\tau\neq-\id$ (with $\pi_2 C(M)$ being non-zero by \cite{Hatcher(1978)}).

\section{Proof of Theorem \ref{thm:main_general_result}}

We may assume that $M$ is connected. Let $H$ be a non-trivial concordance on $S^1\times D^{n-1}$ as given by Theorem \ref{thm:main_result}. Let $i\colon S^1\to M$ be an embedding representing the generator of $\pi_1(M)$. Since $M$ is oriented, $i$ has a trivial normal bundle and induces an embedding $\bar i\colon S^1\times D^{n-1}\to M$. So we may extend $H$ by the identity to a concordance $\bar H$ on $M$, such that $e(\bar H)=\id$.

In the stable range $\dim(M) \geq 9$, the assignment $M\mapsto C(M)$ is a homotopy functor \cite{Hatcher(1978)}. Thus, if $\rho\colon M\to S^1\times D^{n-1}$ classifies (up to homotopy) the universal covering of $M$, then the composite homomorphism
\[C(S^1\times D^{n-1}) \xrightarrow{\bar i_*} C(M) \xrightarrow{\rho_*} C(S^1\times D^{n-1})\]
is the identity. Hence $\bar H=\bar i_*(H)\neq 0$.

\end{document}